\documentclass{amsart}
\usepackage{natbib}
\bibpunct{(}{)}{;}{a}{}{,}

\usepackage[english]{babel}
\usepackage{amssymb, amsthm, amsbsy}
\usepackage{hyperref}

\newcommand{\abbr}[1]{{\sc\lowercase{#1}}}
\theoremstyle{plain}
\newtheorem{thm}{Theorem}[section]
\newtheorem{lem}[thm]{Lemma}
\newtheorem{prop}[thm]{Proposition}

\newtheorem*{thm*}{Theorem}
\theoremstyle{definition}
\newtheorem{conj}{Conjecture}

\newtheorem*{claim*}{Claim}

\newtheorem*{cond*}{}

\newcommand\R{{\mathbb R}}

\newcommand\Z{{\mathbb Z}}

\newcommand\Oh{{\mathcal O}}
\newcommand\oh{{o}}

\newcommand\vep{{\varepsilon}}

\newcommand\Ex{{\mathbb E}}
\newcommand\pr{{\mathbb P}}

\newcommand{\ceil}[1]{\lceil{#1}\rceil}

\newcommand\Bin{\mathop{\sf Bin}\nolimits}

\newcommand{\spec}{\mathsf{Spec}}

\newcommand\gnrn{G({\mathcal X}_n;r(n))}
\newcommand\gdrn{G({\mathcal D}_n;r(n))}

\newcommand\dcb{[0,1]^d}

\newcommand\whp{\text{\bf whp}}

\newcommand\uf{{\mathsf{Unif}}}
\newcommand\nb[1]{{\mathcal{N}}(#1)}
\newcommand\pxn{P(\mathcal{X}_n)}
\newcommand\pdn{P(\mathcal{D}_n)}
\newcommand\spxn{\spec(\mathcal{X}_n)}
\newcommand\spdn{\spec(\mathcal{D}_n)}
\newcommand\lxi{\lambda_i\left(\mathcal{X}_n\right)}
\newcommand\ldi{\lambda_i\left(\mathcal{D}_n\right)}
\newcommand\mat[2]{\mathbb{M}_{#1\times#1}\left(#2\right)}
\newcommand\tr[1]{\mathsf{Tr}\left(#1\right)}
\newcommand\trans[1]{{#1}^{\mathrm{T}}}
\newcommand\hs[1]{\left\|#1\right\|_{\mathrm{HS}}}
\newcommand\ab[2]{\genfrac{}{}{0pt}{}{#1}{#2}}
\newcommand\nbu{\left|\nb{u}\right|}
\newcommand\nui{\left|\nb{u'}\right|}
\newcommand\nuu{\left|\mathcal{N}\left(u,u'\right)\right|}

\newcommand\prob[1]{\mathbb{P}\left\{#1\right\}}
\newcommand\muxn{\mu\left(\mathcal{X}_n\right)}
\newcommand\mudn{\mu\left(\mathcal{D}_n\right)}

\begin{document}
\title[Spectra of random geometric graphs]{The spectrum of a 
random geometric graph is concentrated}
\author[S. Rai]{Sanatan Rai}
\email{sanat@stanford.edu}
\address{Department of Management Science and Engg.,
Stanford University, ca 94305.}

\begin{abstract}
Consider $n$ points distributed uniformly in $[0,1]^d$.  Form a graph
by connecting two points if their mutual distance is no greater than
$r(n)$. This gives a random geometric graph, $\gnrn$, which is
connected for appropriate $r(n)$. We show that the spectral measure of
the transition matrix of the simple random walk (\abbr{srw}) on
$\gnrn$ is concentrated, and in fact converges to that of the graph on
the deterministic grid. 
\end{abstract}
\maketitle

\section{Introduction}
Let $S$ be finite set contained in $\dcb$. Form a graph by connecting
two points, $u,v\in S$ if $\|u-v\|\le r$, obtaining a graph $G(S;r)$.

Let $\mathcal{X}_n$ be a set on $n$ points distributed iid $\uf\dcb$, 
we call $G\sim\gnrn$ a random geomtric graph. The function $r(n)$ is chosen
to be such that $r(n)\downarrow 0$ as $n\uparrow\infty$, but such that
$G$ is a.s. and $\whp$ connected. Herein, $\whp$ denotes with high probability,
ie with a probability greater than $1-n^{-c}$ for some constant $c>0$.
We shall write $\mathcal{D}_n$ for the set of $n$ grid points that are
the intersections of axes parallel lines with  separation $n^{-1/d}$. Hence,
$\gdrn$, is a deterministic graph. 

For a graph $G$, $P(G)$ shall denote the transition probability matrix
for the simple random walk (\abbr{srw}) on $G$. That is:
\[
P(G)_{uv} = \frac{[u\sim v]}{|\nb{u}|},
\]
where $\nb{u}$ is the set of neighbours of vertex $u$, and $[u\sim v]$
is the indicator of the event $u\sim v$. Henceforth, 
we shall write $\pxn$ and $\pdn$ for $P(\gnrn)$ and $P(\gdrn)$,
respectively, with $\spxn$ and $\spdn$ denoting their 
spectra. Furthermore, $\muxn$ and $\mudn$ shall stand for their
spectral measures respectively.

Our main result is:
\begin{thm}\label{thm:ae}
$\spxn$ and $\spdn$ are asymptotically equidistributed $\whp$. Moreover,
for $r(n)$ such that $M_n=\oh(r_n)$ a.s. and $\whp$, 
the random spectral measure $\muxn$:
\[
\muxn(-\infty, x]:=\frac{1}{n}\left|\left\{\lambda\in\spxn: \lambda\le x\right\}\right|
\]
is concentrated, ie, there exists a sequence of deterministic measures 
$\{\mudn\}$, and $c_d>0$ such that ($a(n):=n\pi_dr(n)^d$):
\begin{small}
\begin{gather*}
\pr\left\{\left\|\muxn-\mudn\right\|_{\mathrm{WS}}>\frac{t}{a(n)^{1/4}}\right\}
\le\frac{16na(n)^{1/4}}{t}\times\\
\left[
2\exp\left(-\frac{1}{2}\left(\frac{t^4}{8t^4+4096}\right)^2a(n)\right)
+\exp\left(-c_d\frac{t^8}{512}a(n)\right)\right].
\end{gather*}
\end{small}
Here $M_n$ is the minimum bottleneck matching distance between $\mathcal{X}_n$
and $\mathcal{D}_n$ in $d$-dimensions.
\end{thm}       
The norm $\|\cdot\|_{\mathrm{WS}}$ is the \textbf{Wasserstein distance}:
\[
\left\|\mu-\nu\right\|_{\mathrm{WS}}:=
\sup_{f\text{ Lipschitz}}\left|\int fd\mu -\int fd\nu\right|.
\]
\subsection{Related work}
Eigenvalues of random matrices with iid entries have been studied
extensively, can be said to have begun with the work \cite{wigner},
wherein he showed that the spectral measure of $N^{-1/2}X$ converges
to the semicircle law, whenever $X$ is Hermitian with iid complex
entries. The rate of convergence for this class of random matrices was
computed in \cite{bai}.  The work closest in spirit to ours is the
paper by \cite{GuinZei}, wherein they prove rates of convergence using
concentration of measure ideas.  There is a vast body of research in
this area, and rather than provide a comprehensive review, we only
note here that the results and methods from this literature are not
applicable to our case. In many, an exact expression of the limiting
distributions of eigenvalues is used, and in some the Stiltjes
transform is the main tool. Neither technique seems useful in this
case. On the other hand, our proofs are elementary, and require only a
knowledge of Chernoff-H\"{o}ffding bounds and basic real analysis.

Bounds on the second eigenvalue were obtained by \cite{fmmc}, and used to
establish the mixing times of the \abbr{srw} and the fastest walk on $\gnrn$.

We also note that a Wilandt-Hoffman theorem for the eigenvalues of 
discrete matrix approximations to learning kernels were proved in
\cite{KoltGine, Kolt}.
\section{Proof of the main theorem}
Two sequences of real numbers $\{x_n\}$ and $\{y_n\}$ are said to be
\textbf{asymptotically equidistributed}, if for any $L^1$ function $f$ (\cite{ps}):
\[
\lim_{n\uparrow\infty}\left|
\frac{1}{n}\sum_{i=1}^nf\left(x_i\right)-
\frac{1}{n}\sum_{i=1}^nf\left(y_i\right)\right|=0.
\]

Suppose $\{A_n\}$ and $\{B_n\}$ are two sequencs of matrices, such that
$A_n, B_n\in\mat{n}{\R}$. We say that the two sequences are asymptotically
equivalent iff 
\begin{enumerate}
\item $\|A_n\|$ and $\|B_n\|$ are uniformly bounded, and
\item $\|A_n-B_n\|_{\mathrm{HS}}\rightarrow0$ as $n\uparrow\infty$.
\end{enumerate}
Here $\|\cdot\|$ is the usual operator norm:
\[
\|A\|:=\sup_{\|x\|=1}\left|Ax\right|,
\]
and $\|\cdot\|_{\mathrm{HS}}$ is the \textbf{Hilbert-Schmidt norm}:
\[
\|A\|_{\mathrm{HS}}:=\left(\frac{1}{n}\sum_{i,j}A_{ij}^2\right)^{1/2}=
\left(\frac{1}{n}\tr{\trans{A}A}\right)^{1/2}.
\]
With a few straightforward estimates, our main result shall follow easily 
from (\cite{Gray}):
\begin{thm}\label{thm:graylin}
For any two matrices $A$ and $B$:
\[
\left|\frac{1}{n}\sum_{i=1}^n\lambda_i(A)
-\frac{1}{n}\sum_{i=1}^n\lambda_i(B)\right|\le\hs{A-B}.
\]
Therefore, if $\{A_n\}$ and $\{B_n\}$ are two asymptotically
equivalent sequences of matrices, then their spectra are asymptotically
equivalent. 
\end{thm}
To prove theorem \ref{thm:ae} it suffices to show that $\{\pxn\}$ and 
$\{\pdn\}$ are asymptotically equivalent $\whp$, and that $\hs{\pxn-\pdn}$
is concentrated.

Let $M_n$ denote the length of the minimum bottleneck matching between
$\mathcal{X}_n$ and $\mathcal{D}_n$:
\[
M_n\equiv M_n:=\min_{\ab{\phi:\mathcal{X}_n\rightarrow\mathcal{D}_n}{%
\phi\text{ \rm matching}}}
\max_{u\in\mathcal{X}_n}\left\|u-\phi(u)\right\|.
\] 
Then it is well known that:
\begin{thm}\label{thm:match}
\[M_n=\begin{cases}
\Oh\left(\left(\frac{\log{n}}{n}\right)^{1/d}\right),&\text{ when }d\ge3,
\text{ \rm[\cite{ShorYuk}]}\\
\Oh\left(\left(\frac{\log^{3/2}{n}}{n}\right)^{1/2}\right),&\text{ when }d=2,
\text{ \rm[\cite{LeiShor}]}\\
\Oh\left(\sqrt{\frac{\log{\vep^{-1}}}{n}}\right),&\text{ w.p. }1-\vep,
\text{ when } d=1.\text{ \rm[\cite{GoRaK}]}
\end{cases}
\]
\end{thm}
Henceforth, for $u\in\mathcal{X}_n$ we shall write $u'$ for its
matched point in $\mathcal{D}_n$ under a minimum bottleneck
matching. We shall denote by $\nb{u}$ the set of neighbours of $u$,
and by $\mathcal{N}(u,u')$, the set of neighbours of $u$ that are
mapped to a neighbour of $u'$.  We are now in a position to state a
concentration result:
\begin{lem}\label{lem:hs-conc}
If $r(n)$ is such that $M_n=\oh(r_n)$ a.s. and $\whp$, then there is a constant $c_d>0$:
\begin{small}
\begin{gather*}
\prob{\hs{\pxn-\pdn}^2>\frac{t}{a(n)}}\le2n\times\\
\left[
2\exp\left(-\frac{1}{2}\left(\frac{t}{8t+16}\right)^2a(n)\right)
+\exp\left(-c_d\frac{t^2}{2}a(n)\right)\right]
\end{gather*}
\end{small}
where $a(n):=n\pi_dr(n)^d$.
\end{lem}
\begin{proof}
\begin{eqnarray}
\hs{\pxn-\pdn}^2&=&\frac{1}{n}\sum_{u\in\mathcal{X}_n}\sum_{v\in\mathcal{X}_n}
\left(\frac{[u\sim v]}{\nbu}-\frac{[u'\sim v']}{\nui}
\right)^2\nonumber\\ &=&\frac{1}{n}\sum_{u\in\mathcal{X}_n}\left(
\frac{1}{\nbu}+\frac{1}{\nui}-\frac{2\nuu}{\nbu\nui}
\right)\nonumber\\
&\le&\max_{u}\left(\frac{1}{\nbu}+\frac{1}{\nui}-\frac{2\nuu}{\nbu\nui}\right).\label{eq:nbs}
\end{eqnarray}
So that:
\begin{gather*}
\prob{\hs{\pxn-\pdn}^2>\frac{t}{a(n)}}\le\\
n\max_u\prob{\left(\frac{1}{\nbu}+\frac{1}{\nui}-\frac{2\nuu}{\nbu\nui}\right)
>\frac{t}{a(n)}}
\end{gather*}
Now
\begin{gather*}
\prob{\left(\frac{1}{\nbu}+\frac{1}{\nui}-\frac{2\nuu}{\nbu\nui}\right)
>\frac{t}{a(n)}}\\
\le\prob{\frac{2}{\nui}+\left|\frac{1}{\nbu}-\frac{1}{\nui}\right|
-\frac{2\nuu}{\nbu\nui}>\frac{t}{a(n)}}\\
\le\prob{\left|\frac{1}{\nbu}-\frac{1}{\nui}\right|>\frac{t}{8\nui}}
+\underbrace{\prob{\frac{\nbu-\nuu}{\nbu}>\frac{t}{8}}}_{=:\ddag},
\end{gather*}
since $n\pi_dr(n)^d/4\le\nui\le n\pi_dr(n)^d$. To bound $\ddag$,
observe that if $\|u-v\|\le r(n)-2M_n$, then $\|u'-v'\|\le r(n)$. Thus,
all points within a radius of $r(n)-2M_n$ map to neighbours of $u'$.
Furthermore:
\begin{eqnarray*}
\ddag&\le&\prob{\left|\nbu-\nui\right|>\frac{t}{8t+16}\nui}\\
&+&\prob{\left|\nui-\nuu\right|>\frac{t}{16}\nui}
\end{eqnarray*}
Now $\nuu$ is stochastically greater than $\Bin(n, \pi_dr(n)^d(1-2M_n/n)^d)$,
so that the last term satisfies:
\[
\prob{\left|\nui-\nuu\right|>\frac{t}{16}\nui}\le
2\exp\left(-c_d\frac{t^2}{2}n\pi_dr(n)^d\right),
\]
where $c_d>0$ is such that:
\[
c_d\frac{t}{8}\le\min\left(\frac{1+t/8}{\left(1-2M_n^+/r(n)\right)^d}-1,
\left|\frac{1-t/8}{\left(1-2M_n^+/r(n)\right)^d}-1\right|\right)
\]
for all large enough $n$. Here $M_n^+$ is a constant depending only
on $n$ such that $M_n\le M_n^+$, a.s. and $\whp$. Since we have
chosen $r(n)$ such that $M_n/r(n)\rightarrow0$, such a constant $c_d$ always exists. 

Putting all our estimates together, and using lemma \ref{lem:rec} we obtain:
\begin{small}
\begin{gather*}
\prob{\hs{\pxn-\pdn}^2>\frac{t}{a(n)}}\le2n\times\\
\left[\exp\left(-\frac{1}{2}\left(\frac{t}{t+8}\right)^2a(n)\right)
+\exp\left(-\frac{1}{2}\left(\frac{t}{8t+16}\right)^2a(n)\right)
+\exp\left(-c_d\frac{t^2}{2}a(n)\right)\right]
\end{gather*}
\end{small}
and the conclusion of the lemma follows.
\end{proof}
\begin{prop}\label{prop:ae}
When $r(n)$ is such that $M_n=r(n)$ a.s. and $\whp$, for $g(n)\uparrow\infty$, 
the two sequencs $\{\pxn\}$ and $\{\pdn\}$ are asymptotically
equivalent $\whp$.
\end{prop}
\begin{proof}
Since matrices are stochastic, $\|\pxn\|=\|\pdn\|=1$, and hence they
are of uniformly bounded norm. Lemma \ref{lem:hs-conc} implies that
$\hs{\pxn-\pdn}\rightarrow0$, $\whp$. Hence $\{\pxn\}$ and $\{\pdn\}$
are asymptotically equivalent $\whp$.
\end{proof}
We now prove that the concentration of the Hilbert-Schmidt norm of
the matrices $\pxn$ implies concentration of the spectrum. We first 
need a lemma:
\begin{lem}\label{lem:conc}
\begin{gather*}
\prob{\sup_{f\text{ \rm Lipschitz}}\left|
\frac{1}{n}\sum_{i=1}^nf\left(\lxi\right)
-\frac{1}{n}\sum_{i=1}^nf\left(\ldi\right)\right|>4\vep}\\
\le\frac{2}{\vep}\prob{\hs{\pxn-\pdn}>\vep^2}.
\end{gather*}
\end{lem}
\begin{proof}
When $f$ is linear, the result follows directly from 
theorem \ref{thm:graylin} and lemma \ref{prop:ae}. To obtain a uniform
bound over \emph{all} Lipschitz $f$, we use the trick in \cite{GuinZei} of 
approximating by a finite class of functions, for a given error bound. 

To that end, let $f$ be Lipschitz on $[-1,1]$, and fix $\vep>0$. Define:
$g_\vep(x):=x[0\le x\le\vep]+\vep[x>\vep]$. Then for $f_\vep$
recursively defined by:
\[
f_\vep(x)=\sum_{i=0}^{\ceil{x/\vep}}\left(%
2\left[f\left((i+1)\vep\right)>f\left(i\vep\right)\right]-1\right)%
g_\vep\left(x-i\vep\right),
\]
we must have $\|f-f_\vep\|\le\vep$. Thus we can approximate any
Lipschitz $f$ to within $\vep$ by a weighted sum of at most $2/\vep$
functions, each of which has two constant portions joined by a linear
part. Since the weights are $\pm 1$, we have:
\begin{gather*}
\pr\left\{\sup_{f\text{ \rm Lipschitz}}\left|%
\frac{1}{n}\sum_{i=1}^nf\left(\lxi\right)
-\frac{1}{n}\sum_{i=1}^nf\left(\ldi\right)
\right|>4\vep\right\}\\
\le\frac{2}{\vep}\sup_{k}\pr\left\{%
\left|\frac{1}{n}\sum_{i=1}^ng_k\left(\lxi\right)
-\frac{1}{n}\sum_{i=1}^ng_k\left(\ldi\right)\right|>\vep^2
\right\}\\
\le\frac{2}{\vep}\pr\left\{\hs{\pxn-\pdn}>\vep^2\right\}
\end{gather*}
where $g_k(x):=g_\vep(x-k\vep)$, and the last inequality follows from 
theorem \ref{thm:graylin}.
\end{proof}
\begin{proof}[Proof of therem \ref{thm:ae}]
Proposition \ref{prop:ae} and theorem \ref{thm:graylin} imply that the spectra 
are asymptotically equidistributed $\whp$. By lemma \ref{lem:conc}, and 
the estimate in lemma \ref{lem:hs-conc}:
\begin{small}
\begin{gather*}
\prob{\left\|\muxn-\mudn\right|_{\rm{WS}}>\frac{t}{a(n)^{1/4}}}\\
\le\frac{8a(n)^{1/4}}{t}
\prob{\hs{\pxn-\pdn}>\frac{t^4}{16\sqrt{a(n)}}}\\
\le\frac{16na(n)^{1/4}}{t}\times\\
\left[
2\exp\left(-\frac{1}{2}\left(\frac{t^4}{8t^4+4096}\right)^2a(n)\right)
+\exp\left(-c_d\frac{t^8}{512}a(n)\right)\right]
\end{gather*}
\end{small}
\end{proof}
\section{Discussion and open problems}
Note that if we define $\mu_{n}^{\{x\}}$ to be the empirical measure
\[
\mu_{n}^{\{x\}}(\alpha)=\frac{1}{n}\sum_{i=1}^n\delta(\alpha-x_i),
\]
then $\{x_n\}$ and $\{y_n\}$ are asymptotically equidistributed 
if and only if $\|\mu_n^{\{x\}}-\mu_n^{\{y\}}\|_{\mathrm{WS}}\rightarrow 0$,
since we may approximate integrable functions by Lipschitz ones.

Simulations suggest that in fact, the eigenvalues of $\pxn$ and $\pdn$
are \textbf{asymptotically absolutely equally distributed} (\cite{trench}):
\begin{conj}
For any $f\in\mathcal{C}[-1,1]$:
\[
\lim_{n\uparrow\infty}
\frac{1}{n}\sum_{i=1}^n\left|f\left(\lxi\right)-f\left(\ldi\right)\right|^2
=0,
\]
where $\lxi$ and $\ldi$ are the $i$th largest eigenvalues of $\pxn$  and
$\pdn$ respectively.
\end{conj}
Such a result would have followed immediately, were it known that the matrices
were symmetric. This is of course, not true. However, the matrices are
\emph{almost} symmetric in the asymptotic limit, so one would expect
the matrices to be \emph{almost} diagonalisable via unitary matrices. Then
the conjecture would follow immediately from the theory of asymptotically
equivalent matrices. As it stands, one needs some other manner of bounding
the differences:
\[
\frac{1}{n}\sum_{i=1}^n\left|\lxi-\ldi\right|^2,
\]
to obtain say, a Wilandt-Hoffman type theorem. Then the conjecture would
follow,  the Stone-Weierstrass theorem.

It also seems reasonable to posit exponential tail bounds 
for the Wielandt-Hoffman type result. 

Another natural question to ask is: what is the behaviour of the
resolvent $R_n(x):=(I-x\pxn)^{-1}$, for
$x\not\in\spec(\mathcal{X}_n)$. Whereas, $(I-\pxn)\rightarrow\Delta$,
the asymptotic behaviour of the resolvent is not known. Convergence
results for the resolvents of random walks in random environments,
specifically, on $\Z^d$ with ergodic (random) bond percolation were
proved in \cite{kunn}.  

\appendix\section{Auxilliary lemmata} Let
$X\sim\Bin(n,p)$, then we have the following two tail bounds:
\begin{lem}[Chernoff-H\"{o}ffding bounds for reciprocals]\label{lem:rec}
For $t\ge0$:
\[
\prob{\left|\frac{1}{X}-\frac{1}{\Ex{X}}\right|>\frac{t}{\Ex{X}}}\le2\exp\left(-\frac{1}{2}\left(\frac{t}{1+t}\right)^2\Ex{X}\right).
\]
\end{lem}
\begin{proof}
Note that:
\begin{eqnarray*}
\prob{\left|\frac{1}{X}-\frac{1}{\Ex{X}}\right|>\frac{t}{\Ex{X}}}
&=&\prob{\left|X-\Ex{X}\right|>tX}\\
&\le&\prob{\left|X-\Ex{X}\right|>t\left(\Ex{X}-
\left|X-\Ex{X}\right|\right)}\\
&=&\prob{\left|X-\Ex{X}\right|>\frac{t}{1+t}\Ex{X}}\\
&\le&2\exp\left(-\frac{1}{2}\left(\frac{t}{1+t}\right)^2\Ex{X}\right).
\end{eqnarray*}
\end{proof}
\bibliographystyle{plainnat}
\bibliography{gnrn-ev}
\end{document}